\documentclass[11pt]{amsart}
\usepackage{graphicx}
\usepackage{hyperref}
\usepackage{subfig}

\usepackage{amsthm}
\usepackage[centering,margin=1in]{geometry}

\usepackage[numbers]{natbib}
\usepackage{amsmath,amssymb,enumerate}
\usepackage[all]{xy}
\usepackage{mathtools}
\usepackage{ytableau}
\usepackage{longtable}
\ytableausetup{boxsize=2.4ex,centertableaux}

\theoremstyle{plain} 
\newtheorem{theorem}{Theorem}[section]

\theoremstyle{definition}
\newtheorem{example}[theorem]{Example}
\newtheorem{definition}[theorem]{Definition}

\theoremstyle{remark}
\newtheorem{remark}[theorem]{Remark}

\newlength{\cellsize}
\newcommand\tableau[1]{
\vcenter{
\let\\=\cr
\baselineskip=-16000pt
\lineskiplimit=16000pt
\lineskip=0pt
\halign{&\tableaucell{##}\cr#1\crcr}}}

\newcommand{\tableaucell}[1]{{%
\def \arg{#1}\def \void{}%

\ifx \void \arg
\vbox to \cellsize{\vfil \hrule width \cellsize height 0pt}%
\else
\unitlength=\cellsize
\begin{picture}(1,1)
\put(0,0){\makebox(1,1){$#1$}}
\put(0,0){\line(1,0){1}}
\put(0,1){\line(1,0){1}}
\put(0,0){\line(0,1){1}}
\put(1,0){\line(0,1){1}}
\end{picture}%
\fi}}
\def \g{\mathfrak{g}}

\def\g{{\mathfrak g}}
\def\w{{\rm wt}}

\def \g{\mathfrak{g}}

\begin{document}
\title[Combinatorial representation theory of Lie algebras]{Combinatorial representation theory of Lie algebras.\\ Richard Stanley's work and the way it was continued}
\author{Cristian Lenart} 
\address{Department of Mathematics and Statistics, State University of New York at Albany,
Albany, NY 12222, USA}
\email{clenart@albany.edu}
\thanks{C.L. was partially supported by the NSF grant DMS--1101264, and gratefully acknowledges the hospitality and support of the Max-Planck-Institut f\"ur Mathematik in Bonn, where this article was written.}

\bibliographystyle{plain}

\maketitle

\vspace{-0.5cm}

\begin{center}
{\small \em Dedicated to Richard Stanley on the occasion of his seventieth birthday.}
\end{center}

\vspace{0.3cm}

Representation theory is a fundamental tool for studying group symmetry $-$ geometric, analytic, or algebraic $-$ by means of linear algebra, which has important applications to other areas of mathematics and mathematical physics. One very successful trend in this field in recent decades involves using combinatorial objects to model the representations, which allows one to apply combinatorial methods for studying them, e.g., for concrete computations. This trend led to the emergence of combinatorial representation theory, which has now become a thriving area. 

Richard Stanley played a crucial role, through his work and his students, in the development of this new area. In the early stages, he has the merit to have pointed out to combinatorialists, in \cite{stausa,stagln}, the potential that representation theory has for applications of combinatorial methods. Throughout his distinguished career, he wrote significant articles which touch upon various combinatorial aspects related to representation theory (of Lie algebras, the symmetric group, etc.). I will describe some of Richard's contributions involving Lie algebras, as well as recent developments inspired by them (including some open problems), which attest the lasting impact of his work.

\subsection*{Acknowledgement} 
Meeting Richard was a defining moment of my life and career, as well as of my evolution since then. Therefore, I always said that I consider him my mentor. I admire him both as a person and as a scientist, while 
his vast work in combinatorics and on its many relationships with other areas of mathematics is a continuous  inspiration for me.

\section{Introduction}

The story starts with the well-known fact 
that the {\em Schur symmetric polynomials} $s_\lambda(x_1,\ldots,x_n)$, indexed by partitions $\lambda=(\lambda_1\ge\ldots\ge\lambda_{n}\ge 0)$, are the irreducible polynomial characters of $GL_n=GL_n({\mathbb C})$. (A representation $X\in GL_n\mapsto \psi(X)\in GL(V)$ of $GL_n$ on a vector space $V$ is polynomial, resp. rational, if the entries of $\psi(X)$ are polynomial, resp. rational functions of the entries of $X$.) We have the immediate connection to combinatorics, due to the well-known  formula:
\begin{equation}\label{schurf}
s_\lambda(x_1,\ldots, x_n)=\sum_{T\in SSYT_{[n]}(\lambda)} x^T\,.
\end{equation}
Here $SSYT_{[n]}(\lambda)$ denotes the set of {\em semistandard Young tableaux} of shape $\lambda$ filled with integers in $[n]:=\{1,\ldots,n\}$, that is, fillings of the Young diagram $\lambda$ with weakly increasing rows (from left to right) and strictly increasing columns (from top to bottom), where we use the English convention; the monomial $x^T$ is given by $x_1^{c_1}\ldots x_n^{c_n}$, where $c_i$ is the number of $i$'s in $T$. 

The very rich combinatorics of symmetric polynomials (or symmetric functions, if infinitely many variables are used) is beautifully described by Richard in \cite{staec2}[Chapter 7], to which we refer for more information on the topics below. Many questions in the representation theory of $GL_n$, or the Lie algebra $\mathfrak{sl}_n=\mathfrak{sl}_n({\mathbb C})$, are answered in combinatorial terms. A typical example is the important problem of decomposing tensor products of irreducible representations; this problem is reduced to describing the structure constants (or Littlewood-Richardson coefficients) for the multiplication of Schur functions, and its  solution is given by the celebrated {\em Littlewood-Richardson rule}. 

\section{The Robinson-Schensted-Knuth correspondence}\label{rsk}

\subsection{The classical correspondence} The {\em Robinson-Schensted-Knuth (RSK) correspondence} is a well-known bijection between non-negative integer matrices with finite support and pairs of SSYT of the same shape. It has an important role in the theory of symmetric functions; for instance, it bijectively proves the {\em Cauchy identity} for Schur functions, cf. \eqref{scauchy}. 

Let us briefly explain the RSK correspondence. Given a non-negative integer matrix $(a_{ij})$ with finite support, we first associate with it (in a bijective way) a biword, that is, a pair of words $(\mathbf{i}=i_1\ldots i_r,\: \mathbf{j}=j_1\ldots j_r)$, viewed up to permutation of biletters $(i_k,j_k)$. The bijection is given by considering $a_{ij}$ biletters $(i,j)$ for each $a_{ij}>0$. Unless otherwise specified, we think of a biword as its lexicographic representative, that is, we order the biletters such that $i_k<i_{k+1}$, or $i_k=i_{k+1}$ and $j_k\le j_{k+1}$, for all $k$. We then map $(\mathbf{i},\mathbf{j})$ to a pair $(P,Q)$ of SSYT of the same shape by the {\em Schensted insertion} procedure applied to the word $\mathbf{j}$ (which gives $P$), run in parallel with a ``recording'' procedure applied to the word $\mathbf{i}$ (which gives $Q$). It is well-known that the transpose of the matrix $(a_{ij})$, i.e., the biword $(\mathbf{j},\mathbf{i})$, is mapped to $(Q,P)$. 

\subsection{A skew analogue} In Richard's joint paper with B. Sagan \cite{sasrsa}, a generalization of the RSK correspondence to skew tableaux is constructed. 

 Given a skew Young diagram $\lambda/\mu$, for $\mu\subset\lambda$, we denote its semistandard fillings by $SSYT(\lambda/\mu)$, and define the skew Schur function $s_{\lambda/\mu}(x)$ in infinitely many variables by the analogue of \eqref{schurf}.

\begin{theorem}\label{srsk}{\rm \cite{sasrsa}[Theorem 6.1]} Let $\alpha,\beta$ be fixed partitions. There is a bijection $(({\mathbf i},{\mathbf j}),T,U)\leftrightarrow (P,Q)$ between (lexicographic) biwords $({\mathbf i},{\mathbf j})$ with $T\in SSYT(\alpha/\mu)$, $U\in SSYT(\beta/\mu)$, and $P\in SSYT(\lambda/\beta)$, $Q\in SSYT(\lambda/\alpha)$ such that ${\mathbf i}\cup T=P$ and ${\mathbf j}\cup U=Q$, as unions of the corresponding multisets of entries (defined by adding the number of repetitions of an element).
\end{theorem}

The bijection in Theorem \ref{srsk} is based on two types of insertion procedures: an external and an internal insertion, both of which are based on the usual Schensted insertion. The skew RSK correspondence bijectively proves the following generalization of the Cauchy identity (the classical case corresponds to $\alpha=\beta$ being the empty partition):
\begin{equation}\label{scauchy}
\sum_\lambda s_{\lambda/\beta}(x)\,s_{\lambda/\alpha}(y)=\sum_\mu s_{\alpha/\mu}(x)\,s_{\beta/\mu}(y)\prod_{i,j} (1-x_i y_j)^{-1}\,.
\end{equation}

Other results related to the skew RSK correspondence are derived in \cite{sasrsa}: transposing the matrix exchanges the two output tableaux (like in the classical RSK), various identities for $s_{\lambda/\mu}(x)$ and the number of skew standard Young tableaux $f^{\lambda/\mu}$ (a SSYT is standard if its entries are $1$ through $n=|\lambda/\mu|$), a skew dual RSK correspondence (corresponding to $0,1$-matrices and output tableaux of conjugate shapes), and an RSK correspondence for so-called skew {\em shifted tableaux}. 

\subsection{A connection to representation theory} We will give a representation-theoretic interpretation of the RSK correspondence, in terms of {\em Kashiwara's crystals} \cite{kascbq}. 

Crystals are colored directed graphs encoding representations $V$ of a quantum group $U_q({\mathfrak g})$, i.e., a certain $q$-deformation of the enveloping algebra of a Lie algebra ${\mathfrak g}$, in the limit $q\to 0$. The vertices of the crystal correspond to the elements $b$ of a so-called {\em crystal basis} $B$ of $V$, and the edges give the action on $B$ of the  {\em Kashiwara} (or {\em crystal}) {\em operators} $\widetilde{f}_i$ (i.e., certain analogues of the Chevalley generators $f_i$) in the limit $q\to 0$. More precisely, whenever $B$ exists, $\widetilde{f}_i$ permutes its elements or sends them to $0$, so $\widetilde{f}_i\::\:B\rightarrow B\cup\{\mathbf{0}\}$, and $\widetilde{f}_i(b)=b'$ is represented by an $i$-colored edge $b\stackrel{i}{\longrightarrow}b'$. 

For instance, the $r$th tensor power of the vector representation of the Lie algebra $\mathfrak{gl}_{n}$ has a crystal basis indexed by words of length $r$ in the alphabet $[n]$. By taking the direct limit as $n\to\infty$, we can define the $\mathfrak{gl}_{>0}$-crystal structure on words in the alphabet ${\mathbb Z}_{>0}$. In fact, there is a nice construction of the map $\widetilde{f}_i$ on words (where $i\in {\mathbb Z}_{>0}$), based on the so-called {\em signature rule}. Namely, given a word $w$, form a binary $\pm$-word by recording each $i$ (resp. $i+1$) in $w$ as a $+$ (resp. $-$), from left to right. Then cancel adjacent pairs $-+$ as long as possible, to obtain the $i$-signature of $w$, which has the form $+\ldots+-\ldots-$. If this contains no $+$, then $\widetilde{f}_i(w):=\mathbf{0}$; otherwise, $\widetilde{f}_i(w)$ is defined by changing to $i+1$ the $i$ in $w$ corresponding to the rightmost $+$ in the $i$-signature. For instance, $\widetilde{f}_2(3432122\mathbf{2}342)=3432122\mathbf{3}342$. 

It is useful to identify a SSYT $T$ with its word $w(T)$, by reading the entries of $T$ from left to right in each row, with the rows considered from bottom to top. It is well-known that the class of (words of) SSYT is closed under the crystal operators, so it forms a $\mathfrak{gl}_{>0}$-subcrystal. 

We define (commuting) left and right $\mathfrak{gl}_{>0}$-crystal structures on pairs $(P,Q)$ of SSYT by letting the crystal operators act on $P$, resp. $Q$. We can also define a right $\mathfrak{gl}_{>0}$-crystal structure on biwords $(\mathbf{i},\mathbf{j})$ by ordering the biword lexicographically, as above, which gives $(\mathbf{i}_R,\mathbf{j}_R)$, and by letting the crystal operators act on $\mathbf{j}_R$. Similarly, by exchanging the roles of $\mathbf{i}$ and $\mathbf{j}$ in the definition of the lexicographic order, which gives  $(\mathbf{i}_L,\mathbf{j}_L)$, we can define a left $\mathfrak{gl}_{>0}$-crystal structure on biwords.

\begin{theorem}\label{doublecg}{\rm \cite{lasdcg}} The RSK bijection identifies the left (resp. right) crystal structure on biwords with the right (resp. left) crystal structure on pairs of SSYT. Namely, assuming $(\mathbf{i},\mathbf{j})\leftrightarrow(P,Q)$, we have $\widetilde{f}_i(\mathbf{j}_R)\ne\mathbf{0}$ if and only if $\widetilde{f}_i(P)\ne\mathbf{0}$, in which case
\[(\mathbf{i}_R,\widetilde{f}_i(\mathbf{j}_R))\leftrightarrow(\widetilde{f}_i(P),Q)\,.\]
Similarly, $\widetilde{f}_i(\mathbf{i}_L)\ne\mathbf{0}$ if and only if $\widetilde{f}_i(Q)\ne\mathbf{0}$, in which case
\[(\widetilde{f}_i(\mathbf{i}_L),\mathbf{j}_L)\leftrightarrow(P,\widetilde{f}_i(Q))\,.\]
\end{theorem}

More generally, one can consider the Lie algebra $\mathfrak{gl}_\infty$ of matrices of finite support with rows and columns indexed by ${\mathbb Z}^*$, which has $\mathfrak{gl}_{>0}$ as a subalgebra, together with the obvious isomorphic copy  $\mathfrak{gl}_{<0}$. In this context, Kwon states Theorem \ref{doublecg} in \cite{kwodcg} as an isomorphism of $(\mathfrak{gl}_{<0},\mathfrak{gl}_{>0})$-crystals; here the alphabet corresponding to $\mathfrak{gl}_{<0}$, which indexes the corresponding crystal operators, is taken to be $\mathbb{Z}_{<0}$. 
Moreover, an $\widetilde{f}_0$ crystal operator is defined both for biwords and pairs of SSYT, giving $\mathfrak{gl}_\infty$-crystal structures. Then Theorem \ref{doublecg} is strengthened in \cite{kwodcg} by stating that the RSK correspondence is a $\mathfrak{gl}_\infty$-crystal isomorphism. Note, however, that not the original RSK correspondence is used here, but its seventh  variation (out of eight) in \cite{fulyt}[Appendix A.4.2, Theorem 1].

The RSK bijection is then interpreted in \cite{kwodcg} as an isomorphism between the above $\mathfrak{gl}_\infty$-crystal of matrices and the crystal of a generalized Verma module for $\mathfrak{gl}_\infty$. Finally, a refinement of this isomorphism is given; a special case of it is precisely the Sagan-Stanley skew RSK correspondence. 

A representation-theoretic interpretation of \eqref{scauchy}, also based on crystals, is given in \cite{kwopae}. Furthermore, an application of the skew dual RSK correspondence in \cite{sasrsa}  to the construction of crystal bases of certain modules for the quantum superalgebra $U_q(\mathfrak{gl}(m|n))$ is given in \cite{kwocbq}.

\section{$\mathfrak{sl}_2$-representations and related posets} 

\subsection{Combinatorics from $\mathfrak{sl}_2$-representations}\label{slrep} The representations of the Lie algebra $\mathfrak{sl}_2=\mathfrak{sl}_2(\mathbb{C})$ of complex $2\times 2$ matrices of trace $0$ are a powerful tool for proving various combinatorial properties of certain posets. They have been used by Richard in \cite{stausa,stagln,stalcu}.

We start by observing that the character of an $\mathfrak{sl}_2$-representation on an $n$-dimensional vector space $V$, that is, a Lie algebra homomorphism $\psi\::\:\mathfrak{sl}_2\rightarrow\mathfrak{gl}(V)$, can be expressed as
\[ch(\psi)=\sum_{i=1}^n q^{e_i}\,;\]
here $e_i$ are the eigenvalues of the $n\times n$ matrix $\psi(H)$, with $H=\left(\!\!\begin{array}{cc} 1&0\\0&-1 \end{array}\!\!\right)$. As it is well-known, the irreducible representations of $\mathfrak{sl}_2$ are labeled by $k\ge 0$, so we denote them by $\psi_k$; we also have
\[ch(\psi_k)=q^{-k}+q^{-k+2}+q^{-k+4}+\ldots+q^k\,.\]
This basic fact implies the following main tool for combinatorics.

\begin{theorem}{\rm \cite{stalcu}[Theorem 15]} Let $\psi$ be an $\mathfrak{sl}_2$-representation with character $ch(\psi)=\sum_i b_i q^i$. Then the following two sequences are symmetric and unimodal:
\[\ldots,\,b_{-4},\,b_{-2},\,b_0,\,b_2,\,b_4,\,\ldots\;\;\;\;\;\mbox{and}\;\;\;\;\ldots,\,b_{-3},\,b_{-1},\,b_1,\,b_3,\,\ldots\;.\]
\end{theorem}

Recall that a sequence $a_0,\,a_1,\,\ldots,\,a_d$ is symmetric if $a_i=a_{d-i}$ for all $i$, and unimodal if, for some $i$, we have $a_0\le a_1\le \ldots\le a_i\ge a_{i+1}\ge\ldots\ge a_d$. Several examples were given in \cite{stausa,stawgh,stagln,stalcu}. Typical ones involve posets where the numbers $a_i$ are the rank cardinalities. Interesting examples come from topology, and are intimately related to the {\em hard Lefschetz theorem}; here the $\mathfrak{sl}_2$-representation is on the cohomology or intersection cohomology of an irreducible complex projective variety. Among these examples are the ones involving parabolic quotients of Weyl groups and the generalized lower bound conjecture for simplicial polytopes. These results are surveyed elsewhere in this volume, so we will not address them here. 

Another interesting property of a ranked poset which can be proved via $\mathfrak{sl}_2$-representations is the {\em strong Sperner property}, see \cite{prorsl,stawgh}. This says that the largest subset of the poset containing no $k+1$-element chain has the same size as the union of the $k$ largest ranks, for each $k$ (the classical Sperner property corresponds to $k=1$). As this property is proved in conjunction with the properties discussed above, the following concept is useful: a ranked poset is called {\em Peck} if it is strongly Sperner and its rank sequence is symmetric and unimodal.

Similar results (unimodality, Sperner property) were derived by Richard in \cite{stauls} by using the superalgebra analogue of $\mathfrak{sl}_2$, namely the orthosymplectic superalgebra $\mathfrak{osp}(1,2)$. 

\subsection{Differential posets} As mentioned above, many interesting $\mathfrak{sl}_2$-representations were constructed based on posets. This fact led to the concept of an {\em $\mathfrak{sl}_2$-poset}. To define it, we need some notation. Given a poset $P$ and a field $K$, let $KP$ denote the $K$-vector space with basis $P$. The covering relation in $P$ is denoted by $\lessdot$. If $P$ is ranked, let $P_i$ be the set of elements of rank $i$. 

\begin{definition}{\rm \cite{prorsl}, \cite{stavdp}[Section 4]} {\em A finite ranked poset $P$ of rank $n$ is called an $\mathfrak{sl}_2$-poset if there exist operators $X,Y\::\:KP\rightarrow KP$ satisfying
\begin{equation}\label{defxy} X x=\sum_{y\::\:y\gtrdot x} c_{xy} \,y \;\;\mbox{ for $x\in P$ (where $c_{xy}\in K$)}\,,\;\;\;\;\;   Y x\in KP_{i-1}\;\;\mbox{ for $x\in P_i$}\,,\end{equation}
as well as the following commutation relation in every rank $i$:}
\begin{equation}\label{comm}
(XY-YX)_i=(2i-n)\, Id_i\,.
\end{equation}
\end{definition}

The operator $X$ is called an order-raising operator, while $Y$ is called a lowering operator. Together they define an $\mathfrak{sl}_2$-representation in the usual way, namely $\left(\!\!\begin{array}{cc} 0&1\\0&0 \end{array}\!\!\right)\mapsto X$ and $\left(\!\!\begin{array}{cc} 0&0\\1&0 \end{array}\!\!\right)\mapsto Y$. Using methods similar to those in \cite{stavdp}, Proctor proved the following theorem, cf. Section \ref{slrep}.

\begin{theorem}{\rm \cite{prorsl}[Theorem 1], \cite{stavdp}} A ranked poset $P$ is Peck if and only if it is an $\mathfrak{sl}_2$-poset.
\end{theorem}

An $\mathfrak{sl}_2$-poset is, in fact, one of the variations in \cite{stavdp} of so-called {\em differential posets}. The latter were defined by Richard in the earlier paper  \cite{stadp}, as a generalization of Young's lattice of partitions (under inclusion).   

\begin{definition}{\rm \cite{stadp}} {\em A poset $P$ is called $r$-differential (for some fixed positive integer $r$) if the operators $U,D\::\:KP\rightarrow KP$ given by
\begin{equation}\label{defud}U x:=\sum_{y\::\:y\gtrdot x} y\,,\;\;\;\;\;D x:=\sum_{y\::\:y\lessdot x} y\end{equation}
satisfy the condition}
\begin{equation}\label{comm0}
DU-UD=r\,Id\,.
\end{equation}
\end{definition}

Note the similarity of \eqref{defxy} and \eqref{comm} with \eqref{defud} and \eqref{comm0}, respectively. Many remarkable algebraic and combinatorial properties of differential posets and of their variations are derived by Richard in \cite{stadp,stavdp}. Other researchers continued this study. For instance, Fomin \cite{fomsad} proves a generalization of the RSK correspondence, cf. Section \ref{rsk}, to differential posets (in fact, to his slightly more general class of {\em dual graded graphs}). His related concept of {\em growth diagrams} (see also \cite{staec2}[Appendix 1]) has proved to be very useful in representation theory in relation to crystals, see \cite{lenccg2,leuajt}. 

Degree 3 relations between the down and up operators on posets defined in \eqref{defud} were considered in the literature. Benkart and Roby \cite{bardua} study the infinite-dimensional associative algebras generated by such operators, which they call {\em down-up algebras}. They show that these algebras exhibit many of the important features of the universal enveloping algebra of $\mathfrak{sl}_2$, including a  Poincar\'e-Birkhoff-Witt type basis and a well-behaved representation theory. 

\subsection{Representation diagrams} We will now discuss the way in which several $\mathfrak{sl}_2$-posets can be pieced together, in order to describe a representation of any semisimple Lie algebra $\mathfrak{g}$, with {\em Chevalley generators} $\{X_i,\,Y_i,\,H_i\}_{i\in I}$. These ideas are due to R. Proctor and his student, R. Donnelly, and we refer to \cite{donepb} for more details. 

We start with a finite ranked poset $P$, viewed as a directed graph with edges $x\rightarrow y$ for each covering relation 
$x\lessdot y$. These edges will be colored by the set $I$ above, so we denote them by $x\stackrel{i}{\rightarrow}y$ for some $i\in I$.  The connected components of the subgraph with edges colored $i$ are called {\em $i$-components}. Beside a given color, each edge $x\rightarrow y$ is labeled with two complex coefficients, which are not both $0$, and which are denoted by $c_{yx}$ and $d_{xy}$. 

Based on the above structure, we define operators $X_i$ and $Y_i$ on $\mathbb{C}P$ for $i$ in $I$, as follows:
\begin{equation}\label{updown}X_i\, x:=\sum_{y\::\:x\stackrel{i}{\rightarrow}y}c_{yx}\,y\,,\;\;\;\;\;\;\;\;\;\;\;Y_i\, y:=\sum_{x\::\:x\stackrel{i}{\rightarrow}y}d_{xy}\,x\,.\end{equation}
For each vertex $x$ of $P$, we also define a set of integers $\{m_i(x)\}_{i\in I}$ by $m_i(x):=2\rho_i(x)-l_i(x)$, where $l_i(x)$ is the rank of the $i$-component containing $x$, and $\rho_i(x)$ is the rank of $x$ within that component. 

The goal is to define a representation of $\mathfrak g$ on $\mathbb{C}P$ by letting the Chevalley generators $X_i$ and $Y_i$ act as in (\ref{updown}), and by setting
\begin{equation}\label{acth}H_i\, x:=m_i(x)\, x\,.\end{equation}
Note the similarity of \eqref{updown} and \eqref{acth} with \eqref{defxy} and \eqref{comm}, respectively (recall that $H_i=[X_i,Y_i]$). We will give necessary and sufficient local conditions on the edge labels, which are translations of the {\em Serre relations} for the Chevalley generators of $\g$.

Let $\{\omega_i\}_{i\in I}$ and $\{\alpha_i\}_{i\in I}$ denote the {fundamental weights} and {simple roots} of the root system corresponding to $\g$, respectively. 
We assign a weight to each vertex of $P$ by $\w(x):=\sum_{i\in I}m_i(x)\omega_i$. We say that the edge-colored poset $P$ satisfies the {\em structure condition} for $\mathfrak{g}$ if 
\begin{equation}\label{strcond}\w(x)+\alpha_i=\w(y)\;\;\;\;\mbox{ whenever $x\stackrel{i}{\rightarrow}y$}\,.\end{equation} 

We now define two conditions on the pairs of edge labels $(c_{yx},d_{xy})$. We call $\pi_{xy}:=c_{yx}\,d_{xy}$ an {edge product}. We say that the edge-labeled poset $P$ satisfies the {\em crossing condition} if for any vertex $y$ and any color $i$ we have
\begin{equation}\label{crossing}\sum_{x\::\:x\stackrel{i}{\rightarrow}y}\pi_{xy}-\sum_{z\::\:y\stackrel{i}{\rightarrow}z}\pi_{yz}=m_i(y)\,.\end{equation}
We say that the edge-labeled poset $P$ satisfies the {\em diamond condition} if for any pair of vertices $(x,y)$ of the same rank and any pair of colors $(i,j)$, possibly $i=j$, we have
\begin{equation}\label{diamond}\sum_{z\::\:x\stackrel{j}{\rightarrow}z\;{\rm and}\;y\stackrel{i}{\rightarrow}z}c_{zx}\,d_{yz}=\sum_{u\::\:u\stackrel{i}{\rightarrow}x\;{\rm and}\;u\stackrel{j}{\rightarrow}y}d_{ux}\,c_{yu}\,.\end{equation}

\begin{theorem}\label{cond}{\rm \cite[Proposition 3.4]{donepb}}
Given an edge-colored and edge-labeled ranked poset $P$, the actions {\rm (\ref{updown})} and {\rm (\ref{acth})} define a representation of $\mathfrak g$ on $\mathbb{C}P$ if and only if $P$ satisfies the structure, crossing, and diamond conditions, namely {\rm \eqref{strcond}}, {\rm \eqref{crossing}}, and {\rm \eqref{diamond}}. 
\end{theorem}

If the given conditions hold, we obtain a {\em weight basis} of the given representation, and the edge-colored poset $P$ together with its edge-labels is called a {\em representation diagram}. 

Based on representation diagrams, in \cite{halccw} we give a short proof of the {\em Gelfand-Tsetlin basis} construction for the irreducible representations of $\mathfrak{sl}_n$; recall that this basis is the only one with respect to which the representation matrices of the Chevalley generators are given by explicit formulas. 
The relevant poset is the set of SSYT of a fixed shape with entries in $[n]$, where an edge colored $i$ corresponds to changing a single entry $i+1$ to $i$ in a SSYT; this is known as the Gelfand-Tsetlin (distributive) lattice. The edge coefficients are certain rational numbers expressed in a nice factored form in terms of the corresponding SSYT. Based on Theorem \ref{cond}, our proof consists of verifying that these coefficients satisfy  {\rm \eqref{crossing}} and {\rm \eqref{diamond}}.
 
Similar explicit constructions of certain symplectic and orthogonal irreducible representations (indexed by rows or columns) were given in \cite{donecf,dlpcro}, as well as of the adjoint representations of all simple Lie algebras \cite{doneba}. Several properties of such bases, such as minimality with respect to the edge set, were also exhibited \cite{donepb}. However, much more work remains to be done, for instance to extend the mentioned results to an arbitrary symplectic or orthogonal irreducible representation. 

\section{The stable behavior of some graded multiplicities} This section is based on Richard's paper \cite{stasbs}, which studies the decomposition into irreducible components of polynomial functions on the adjoint representation of $G=GL_n(\mathbb{C})$, i.e., the representation on its Lie algebra $\mathfrak{g}=\mathfrak{gl}_n(\mathbb{C})$. (In fact, the setup of the paper is more general, also including a reference to the {\em $q$-Dyson conjecture}, and $SL_n(\mathbb{C})$ is used instead.)

The adjoint representation $\mbox{ad}:G\rightarrow GL(\mathfrak{g})$ is the rational representation defined by $\mbox{ad}(X)\,A=XAX^{-1}$. This extends to a representation on the symmetric algebra $S(\mathfrak{g})$, viewed as a graded algebra in the usual way. By a theorem of Kostant, $S(\mathfrak{g})$ is a free module over the ring of $G$-invariant polynomials, with  $S(\mathfrak{g})\cong S(\mathfrak{g})^G\otimes \mathcal{H}_{\mathfrak{g}}$, where $\mathcal{H}_{\mathfrak{g}}$ are the {\em $G$-harmonic polynomials}. 

Recall (see, e.g., \cite{stertt}) that the irreducible rational representations of $G$ are labeled by weakly decreasing sequences of integers with length $n$ or, equivalently, by pairs of partitions $\alpha=(\alpha_1\ge\ldots\ge\alpha_k>0)$, $\beta=(\beta_1\ge\ldots\ge\beta_l>0)$ with $k+l\le n$; the corresponding bijection is given by 
\[(\alpha,\beta)\,\mapsto\,(\alpha_1,\ldots,\alpha_k,0,\ldots,0,-\beta_l,\ldots,-\beta_1)=:[\alpha,\beta]_n\,.\]
We may decompose $\mathcal{H}_{\mathfrak{g}}$ into irreducible components such that each one is homogeneous. For a fixed pair $(\alpha,\beta)$, let $d_i$ be the degrees (also called {\em generalized exponents}) of the components isomorphic to the irreducible representation labeled by $[\alpha,\beta]_n$, and define
\[G_{[\alpha,\beta]_n}(q):=\sum_{i} q^{d_i}\,.\]

Following an idea of Ranee Gupta, who showed that 
\[G_{\alpha\beta}(q):=\lim_{n\to\infty} G_{[\alpha,\beta]_n}(q)\]
exists as a formal power series, Richard expressed the stable graded multiplicities based on the internal product $*$ of Schur functions. This corresponds to the tensor product of irreducible representations of the symmetric group (via the {\em Frobenius characteristic}). Here we define the internal product formally, in the basis of power sum symmetric functions $p_\lambda:=p_{\lambda_1} p_{\lambda_2}\ldots=p_1^{m_1}p_2^{m_2}\ldots$, by
\[p_\lambda*p_\mu=\delta_{\lambda\mu}\, z_\lambda\, p_\lambda\,,\;\;\;\;\;\mbox{where $z_\lambda:=(m_1!\,1^{m_1})(m_2!\,2^{m_2})\ldots$}\,.\]

\begin{theorem}{\rm \cite{stasbs}[Proposition 8.1]}\label{stabk} We have
\[G_{\alpha\beta}(q)=s_\alpha*s_\beta(q,q^2,\ldots)\,.\]
\end{theorem}

Stembridge \cite{stertt} developed the combinatorics of the rational representations of $GL_n$, and used it to express the decomposition of the $k$th tensor power of $\mathfrak{gl}_n$ into modules which are simultaneously irreducible with respect to $GL_n$ and the symmetric group $S_k$. Based on this, he was able to rederive some of Richard's results in \cite{stasbs}.

Howe, Tan, and Willenbring \cite{htwsgm} showed that Theorem \ref{stabk} is a special case of a general stability result, which we now briefly describe. Let $G$ be a complex reductive algebraic group and $K\subset G$ the fixed-point set of a regular involution on $G$. Denote by $\mathfrak{k}\subset\mathfrak{g}$ the corresponding Lie algebras, and by $\mathfrak{p}$ the Cartan complement of $\mathfrak{k}$ in $\mathfrak{g}$. The adjoint action of $K$ on $\mathfrak{p}$ extends to a $K$-module structure on the symmetric algebra $S(\mathfrak{p})$. By the Kostant-Rallis theorem (which generalizes Kostant's theorem mentioned above), we have $S(\mathfrak{p})\cong S(\mathfrak{p})^K\otimes \mathcal{H}_{\mathfrak{p}}$. The main result of \cite{htwsgm} describes the graded multiplicity of the $K$-module $\mathcal{H}_{\mathfrak{p}}$ in the case where $(G,K)$ is a {\em classical symmetric pair} and the degree $d$ is restricted to the so-called {\em stable range}. Theorem \ref{stabk} corresponds to $K=GL_n$ diagonally embedded into $G=GL_n\times GL_n$, and the stable range is $n\ge 2d$. There are nine more classical symmetric pairs, such as $(GL_{n+m},\,GL_n\times GL_m)$, $(SO_{2n},\,GL_n)$,  and $(Sp_{2n},\,GL_n)$, for which the graded multiplicities in the stable range are expressed in terms of the Littlewood-Richardson coefficients. 

\section{Jack, Macdonald, and Hall-Littlewood symmetric polynomials}

In this section I will refer to Richard's paper \cite{stascp} on {\em Jack symmetric functions}, as well as to some of the many related developments; indeed, this paper has an impressive number of citations. 

\subsection{Richard's work on Jack symmetric functions} The Jack functions $J_\lambda^{(\alpha)}(x)$, where $\lambda$ is a partition, are a remarkable basis of symmetric functions depending on a parameter $\alpha$. They can be defined as the eigenfunctions of the {\em generalized Laplace-Beltrami operator}, and for $\alpha=1$ they specialize (up to normalization) to the Schur functions \cite{macsfh}. 

Richard's paper \cite{stascp} is the first one devoted to extending the properties of Schur functions to Jack functions. Some of its results involve: specializations of Jack functions, formulas for special Jack functions and special coefficients, a Cauchy identity (cf. \eqref{scauchy}), skew Jack functions and a related coproduct formula, a {\em Pieri formula} (for the expansion in the basis of Jack functions of the product $J_\lambda^{(\alpha)}(x) J_{(n)}^{(\alpha)}(x)$ by a ``one-rowed'' Jack function), and a formula for the skew Jack functions in terms of skew SSYT (which is derived from the Pieri formula and generalizes \eqref{schurf}). Richard also states an important conjecture: after suitable renormalization, the structure constants for Jack functions (or analogues of Littlewood-Richardson coefficients for Schur functions) belong to $\mathbb{N}[\alpha]$. 

\subsection{Recent related developments} In \cite{kasrcf}, Knop and Sahi prove that the Jack structure constants are in $\mathbb{Z}[\alpha]$, and give a monomial formula for $J_\lambda^{(\alpha)}(x)$ in terms of more general fillings of the Young diagram $\lambda$ than SSYT (see Section \ref{hlsec}). The coefficients in this formula are in $\mathbb{N}[\alpha]$, whereas in Richard's formula they are rational functions in $\alpha$. 

In \cite{cajalb}, by studying the Laplace-Beltrami operator, Cai and Jing give another monomial formula (with rational function coefficients) for $J_\lambda^{(\alpha)}(x)$, as well as a (non-positive) formula for the Jack structure constants. In another paper, Cai and Jing \cite{cajjvo} prove a special case of Richard's conjecture by deriving the generalization of the Pieri rule in which the one-rowed partition is replaced by a rectangle (or a rectangle and another row). Based on this formula, they settle the long-standing problem of realizing any $J_\lambda^{(\alpha)}(x)$ using {\em vertex operators}. 

In \cite{gahmph}, Graham and Hunziker give an interesting application of Richard's conjecture to relating tensor products $V_\lambda\otimes V_\mu$ of irreducible representations of a complex reductive algebraic group $K$ to the product of isotypical components $S_\lambda\cdot S_\mu$ in the representation of $K$ on the symmetric algebra of some representation $V$. In the context of {\em Hermitian symmetric pairs}, they show that, if $V_\nu$ appears in $V_\lambda\otimes V_\mu$, then $S_\nu\subseteq S_\lambda\cdot S_\mu$ in all the classical cases, provided that Richard's conjecture is true.

\subsection{Macdonald polynomials}\label{smacd} Jack polynomials are specializations of {\em Macdonald polynomials} $P_\lambda(x_1,\ldots,x_n;q,t)$, which are symmetric polynomials depending on two parameters $q,t$ \cite{macsfh}; namely, the former are retrieved by setting $q=t^\alpha$ and letting $t\to 1$. In fact, the mentioned Macdonald polynomials correspond to the root system of type $A_{n-1}$, whereas there are versions of them for an arbitrary finite root system. Macdonald polynomials are also eigenfunctions of certain operators.

A breakthrough came with two monomial formulas for Macdonald polynomials: the {\em Haglund-Haiman-Loehr (HHL) formula} \cite{hhlcfm} in type $A$, and the {\em Ram-Yip formula} \cite{raycfm} in arbitrary type. The former is in terms of certain fillings of the Young diagram $\lambda$ and is shorter, while the latter is in terms of extensions of the {\em alcove model} \cite{lapawg,lapcmc}, or {\em LS-gallery model} \cite{gallsg}, in the representation theory of Lie algebras. The setup of the Ram-Yip formula allows one to also derive a very general, albeit involved, Littlewood-Richardson multiplication formula \cite{yiplrr}. Moreover, it was shown in \cite{yiplrr} that, in type $A$, the Pieri special case of this formula ``compresses'' to Macdonald's more efficient Pieri formula \cite{macsfh}. This raises the possibility of deriving, in a similar way, other efficient multiplication formulas for Macdonald polynomials and their specializations (such as Jack polynomials).

\subsection{Hall-Littlewood polynomials}\label{hlsec} We briefly explain the compression phenomenon mentioned above and the related setup by referring to the $q=0$ specialization of Macdonald polynomials, namely to the {\em Hall-Littlewood polynomials}. 
The reason we chose this specialization is twofold: more is known about compression, and the description is easier. Nevertheless, we can draw a parallel to Jack polynomials, while a compression result for Macdonald polynomials is given in \cite{lencfm}.

We start by explaining the corresponding Ram-Yip formula, which was, in fact, given previously in \cite{ramawh,schghl}. As we only refer to type $A_{n-1}$, we will not use a general root system terminology. 

Fix a partition $\lambda=(\lambda_1\ge\ldots\ge\lambda_{n}=0)$. The {\em Bruhat order} $\le$ on permutations in $S_n$ is defined by its covers: $u\lessdot w=u\,(i,j)$ if $\ell(w)=\ell(u)+1$, where $\ell(\,\cdot\,)$ is the length and $(i,j)$ is a transposition. Let $S_n^\lambda:=\{w\in S_n : \mbox{$w(i)<w(i+1)$ for $\lambda_i=\lambda_{i+1}$}\}$. Define the sequence of transpositions
\[\Gamma_k:=\Gamma_{k,k}\Gamma_{k,k-1}\ldots\Gamma_{k,1}\,,\;\;\;\;\mbox{where }\:\Gamma_{k,i}=((i,k+1),\,(i,k+2),\,\ldots,\,(i,n))\,.\]
Define the sequence $\Gamma_\lambda:=\Gamma_{\lambda_1'}\Gamma_{\lambda_2'}\ldots$, called a {\em $\lambda$-chain}, where $\lambda'$ is the conjugate of $\lambda$. The objects of the alcove model are subwords $T=T_1T_2\ldots$ of $\Gamma_\lambda$, where $T_i$ is a subword of $\Gamma_{\lambda_i'}$. More precisely, we consider the collection of {\em admissible pairs} $(w,T)$, with $w\in S_n$ and $T=(t_1,\,\ldots,\,t_m)$ as before:
\begin{equation}\label{defadm}{\mathcal A}(\lambda):=\{(w,T):\mbox{$wt_1\ldots t_{i-1}<wt_1\ldots t_i$ for all $i$, and $wT\in S_n^\lambda$}\}\,,\;\mbox{where $wT:=wt_1\ldots t_m$}.\end{equation}

Define the permutations $\pi_j=\pi_j(w,T):=wT_1\ldots T_j$, where the notation means right multiplication by the transpositions in $T_1,\,\ldots,\, T_j$, cf. \eqref{defadm}. With $(w,T)$ in ${\mathcal A}(\lambda)$ we associate the filling $\sigma(w,T)$ of $\lambda$ whose $j$th column is filled with $\pi_j(1),\,\ldots,\,\pi_j(\lambda_i')$ from top to bottom. Clearly, the rows of $\sigma(w,T)$ are weakly increasing, but the columns are not necessarily increasing. 

\begin{example}\label{ex1} Let $n=4$ and $\lambda=(2,1,0,0)$, so $\lambda'=(2,1)$. We have 
\begin{equation}\label{exlam}\Gamma_\lambda=((2,3),\,(2,4),\,\underline{(1,3)},\,(1,4)\,|\,(1,2),\,(1,3),\,\underline{(1,4)})\,,\end{equation}
where $|$ indicates the concatenation of various $\Gamma_{\lambda_i'}$. Let $T=T_1T_2$ with $T_1=((1,3))$ and $T_2=((1,4))$, cf. the underlines in \eqref{exlam}. Considering $w$ to be the identity in $S_4$, we note that $(w,T)\in{\mathcal A}(\lambda)$; indeed, the conditions in \eqref{defadm} are verified (we use the one-line notation for permutations):
\[w=1234\,<\,w(1,3)=3214=\pi_1\,<\,w(1,3)(1,4)=wT=4213=\pi_2\in S_4^\lambda\,.\]
Hence, the associated filling $\sigma(w,T)$ and the corresponding monomial $x^{\sigma(w,T)}$, cf. \eqref{schurf}, are
\[\sigma(w,T)=\tableau{{3}&{4}\\{2}}\,,\;\;\;\;\;x^{\sigma(w,T)}=x_2x_3x_4\,.\]
\end{example}

\begin{theorem}{\rm \cite{ramawh,schghl}} We have
\begin{equation}\label{schw} P_\lambda(x_1,\ldots,x_n;t)=\sum_{(w,T)\in{\mathcal A}(\lambda)} t^{\frac{1}{2}(\ell(w)+\ell(wT)-|T|)}\,(1-t)^{|T|}\,x^{\sigma(w,T)}\,.\end{equation}
\end{theorem}

\begin{remark} By setting $t=0$ in \eqref{schw}, the only terms that survive correspond to $(w,T)$ for which $w$ is the identity and the chain in Bruhat order in \eqref{defadm} is saturated (i.e., all the steps are covers). The restriction of the map $(w,T)\mapsto \sigma(w,T)$ to these admissible pairs is a bijection to $SSYT_{[n]}(\lambda)$. As it is well-known that $P_\lambda(x_1,\ldots,x_n;t=0)$ is the Schur polynomial $s_\lambda(x_1,\ldots,x_n)$, we recover \eqref{schurf}. Hence the alcove model gives a new way to realize SSYT as saturated chains in Bruhat order.
\end{remark} 

We will now state a formula which is essentially the specialization of the HHL formula to $q=0$. In this sense, it is the analogue of the Knop-Sahi formula for Jack polynomials mentioned above, which is also obtained by specializing the HHL formula, see \cite{hhlcfm}[Section 8]. 

Let ${\mathcal T}_{[n]}(\lambda)$ denote the fillings $\sigma$ of $\lambda$ with entries in $[n]$ satisfying the following conditions: (i)~the rows of $\sigma$ are weakly increasing; (ii) $\sigma(u)\ne\sigma(v)$ whenever the boxes $u$ and $v$ are in the same column or in two consecutive columns, with the right box strictly below the left one (such positions are called {\em attacking}); (iii) the segment of the leftmost column corresponding to rows $i,i+1,\ldots, j$ with $\lambda_i=\lambda_{i+1}=\ldots=\lambda_j$ is (strictly) increasing. 

\begin{theorem}{\rm \cite{hhlcfm,lenhlp}} We have
\begin{equation}\label{hhlmod} P_{\lambda}(x_1,\ldots,x_n;t)
    = \sum_{\sigma\in {\mathcal T}_{[n]}(\lambda) }
    t^{{\rm cinv} (\sigma )}\,(1-t)^{{\rm des}(\sigma)}\: x^{\sigma }\,,
\end{equation}
where ${\rm cinv}(\,\cdot\,)$ and ${\rm des}(\,\cdot\,)$ are certain statistics on fillings (see {\rm \cite{hhlcfm,lenhlp}}). 
 \end{theorem}
 
There is also an older formula for (skew) Hall-Littlewood polynomials, due to Macdonald \cite{macsfh}, which is in terms of (skew) SSYT. This is the analogue of Richard's formula for (skew) Jack polynomials mentioned above. Recall that a tableau $\sigma\in SSYT_{[n]}(\lambda)$ can be represented by a chain of partitions $\lambda^{(0)}\subseteq\lambda^{(1)}\subseteq\ldots\subseteq\lambda^{(n)}=\lambda$, where $\lambda^{(i)}$ is the shape of the subtableau consisting of entries $\le i$. Let $m_i(\lambda):=\lambda_i'-\lambda_{i+1}'$, where $\lambda'$ is the conjugate partition, as usual.

\begin{theorem}{\rm \cite{macsfh}} We have
\begin{equation}\label{machl} P_{\lambda}(x_1,\ldots,x_n;t)=\sum_{\sigma \in SSYT_{[n]}(\lambda)}  \frac{\varphi_\sigma(t)}{b_{\lambda}(t)}\,x^{\sigma}\,,\end{equation}
where
$$ b_{\lambda}(t):= \prod_{i \geq 1}  \varphi_{m_i(\lambda)}(t)\,,\;\;\;\;\mbox{with }\;\varphi_k(t):= (1-t) \ldots (1-t^k)\,, \;\;\mbox{and}$$
$$ \varphi_{\sigma}(t):= \prod_{i=1}^{n} \varphi_{\lambda^{(i)}/\lambda^{(i-1)}}(t),\;\;\mbox{with}\;\, \varphi_{\lambda/\mu}(t):= \prod_{j \in I}1-t^{m_j(\mu)}\;\,\mbox{and}\;\, I:= \{j : \lambda_j'-\mu_j'  > \lambda_{j+1}'-\mu_{j+1}'\}.$$
\end{theorem} 

 We state a first compression result for Hall-Littlewood polynomials. For simplicity, from now on we assume that $\lambda$ has no repeated parts, so $S_n^\lambda=S_n$; for the general case, see \cite{lenhlp}[Section 6].
 
 \begin{theorem}{\rm \cite{lenhlp}}\label{comp1} The map $(w,T)\!\mapsto\! \sigma(w,T)$ is a surjection from ${\mathcal A}(\lambda)$ to ${\mathcal T}_{[n]}(\lambda)$. For all $\sigma$ in ${\mathcal T}_{[n]}(\lambda)$, the sum in {\rm \eqref{schw}} restricted to $(w,T)$ with $\sigma(w,T)\!=\!\sigma$ equals the corresponding term in {\rm \eqref{hhlmod}}. 
\end{theorem}

We now present a stronger compression for Hall-Littlewood polynomials. In fact, this factors through the first compression, via the surjections ${\mathcal A}(\lambda)\rightarrow {\mathcal T}_{[n]}(\lambda)\rightarrow SSYT_{[n]}(\lambda)$ given by $(w,T)\mapsto \sigma(w,T)\mapsto \overline{\sigma}(w,T)$, where $\overline{\sigma}(w,T)$ is obtained from $\sigma(w,T)$ by sorting columns, cf.  Example \ref{ex1}.

\begin{theorem}{\rm \cite{klogmf}}\label{comp2} 
For all $\sigma$ in ${SSYT}_{[n]}(\lambda)$, the sum in {\rm \eqref{schw}} restricted to $(w,T)$ with $\overline{\sigma}(w,T)\!=\!\sigma$ equals the corresponding term in {\rm \eqref{machl}}. 
\end{theorem}

As stated in Section \ref{smacd}, it would be very interesting to derive compressed formulas similar to those in this section for the multiplication structure constants corresponding to Macdonald polynomials or certain specializations of them (e.g., Hall-Littlewood or Jack polynomials). The starting point would be the general Littlewood-Richardson rule in \cite{yiplrr} in terms of the alcove model. Generalizing some of these results to other Lie types would also be very useful, cf. \cite{klogmf,lenhhl}.


\vspace{-1mm}

{\small

\begin{thebibliography}{10}

\bibitem{bardua}
G.~Benkart and T.~Roby.
\newblock Down-up algebras.
\newblock {\em J. Algebra}, 209:305--344, 1998.

\bibitem{cajalb}
W.~Cai and N.~Jing.
\newblock Applications of a {L}aplace-{B}eltrami operator for {J}ack
  polynomials.
\newblock {\em European J. Combin.}, 33:556--571, 2012.

\bibitem{cajjvo}
W.~Cai and N.~Jing.
\newblock Jack vertex operators and realization of {J}ack functions.
\newblock {\em J. Algebraic Combin.}, 39:53--74, 2014.

\bibitem{donecf}
R.~Donnelly.
\newblock Explicit constructions of the fundamental representations of the
  symplectic {L}ie algebras.
\newblock {\em J. Algebra}, 233:37--64, 2000.

\bibitem{donepb}
R.~Donnelly.
\newblock Extremal properties of bases for representations of semisimple {L}ie
  algebras.
\newblock {\em J. Algebraic Combin.}, 17:255--282, 2003.

\bibitem{doneba}
R.~Donnelly.
\newblock Extremal bases for the adjoint representations of the simple {L}ie
  algebras.
\newblock {\em Comm. Algebra}, 34:3705--3742, 2006.

\bibitem{dlpcro}
R.~Donnelly, S.~Lewis, and R.~Pervine.
\newblock Constructions of representations of {$\mathfrak{o}(2n+1,\ \mathbb
  C)$} that imply {M}olev and {R}einer-{S}tanton lattices are strongly
  {S}perner.
\newblock {\em Discrete Math.}, 263:61--79, 2003.

\bibitem{fomsad}
S.~Fomin.
\newblock Schensted algorithms for dual graded graphs.
\newblock {\em J. Algebraic Combin.}, 4:5--45, 1995.

\bibitem{fulyt}
W.~Fulton.
\newblock {\em Young Tableaux}, volume~35 of {\em London Math. Soc. Student
  Texts}.
\newblock Cambridge Univ. Press, Cambridge and New York, 1997.

\bibitem{gallsg}
S.~Gaussent and P.~Littelmann.
\newblock {LS}-galleries, the path model and {MV}-cycles.
\newblock {\em Duke Math. J.}, 127:35--88, 2005.

\bibitem{gahmph}
W.~Graham and M.~Hunziker.
\newblock Multiplication of polynomials on {H}ermitian symmetric spaces and
  {L}ittlewood-{R}ichardson coefficients.
\newblock {\em Canad. J. Math.}, 61:351--372, 2009.

\bibitem{hhlcfm}
J.~Haglund, M.~Haiman, and N.~Loehr.
\newblock A combinatorial formula for {M}acdonald polynomials.
\newblock {\em J. Amer. Math. Soc.}, 18:735--761, 2005.

\bibitem{halccw}
P.~Hersh and C.~Lenart.
\newblock Combinatorial constructions of weight bases. {T}he
  {G}elfand-{T}setlin basis.
\newblock {\em Electron. J. Combin.}, 17:\#R33, 2010.

\bibitem{htwsgm}
R.~Howe, E.-C. Tan, and J.~Willenbring.
\newblock The stability of graded multiplicity in the setting of the
  {K}ostant-{R}allis theorem.
\newblock {\em Transform. Groups}, 13:617--636, 2008.

\bibitem{kascbq}
M.~Kashiwara.
\newblock On crystal bases of the {$q$}-analogue of universal enveloping
  algebras.
\newblock {\em Duke Math. J.}, 63:465--516, 1991.

\bibitem{klogmf}
I.~Klostermann.
\newblock Generalization of the {M}acdonald formula for {H}all-{L}ittlewood
  polynomials.
\newblock {\em Adv. Math.}, 248:1366--1403, 2013.

\bibitem{kasrcf}
F.~Knop and S.~Sahi.
\newblock A recursion and a combinatorial formula for {J}ack polynomials.
\newblock {\em Invent. Math.}, 128:9--22, 1997.

\bibitem{kwodcg}
J.-H. Kwon.
\newblock Demazure crystals of generalized {V}erma modules and a flagged {RSK}
  correspondence.
\newblock {\em J. Algebra}, 322:2150--2179, 2009.

\bibitem{kwopae}
J.-H. Kwon.
\newblock A plactic algebra of extremal weight crystals and the {C}auchy
  identity for {S}chur operators.
\newblock {\em J. Algebraic Combin.}, 34:427--449, 2011.

\bibitem{kwocbq}
J.-H. Kwon.
\newblock Crystal bases of {$q$}-deformed {K}ac modules over the quantum
  superalgebras {$U_q(\mathfrak{gl}(m\vert n))$}.
\newblock {\em Int. Math. Res. Not.}, (2):512--550, 2014.

\bibitem{lasdcg}
A.~Lascoux.
\newblock Double crystal graphs.
\newblock In {\em Studies in memory of Issai Schur (Chevaleret/Rehovot, 2000)},
  volume 210 of {\em Progr. Math.}, pages 95--114. Birkh\"auser Boston, Boston,
  MA, 2003.

\bibitem{lenccg2}
C.~Lenart.
\newblock On the combinatorics of crystal graphs. {II}. {T}he crystal commutor.
\newblock {\em Proc. Amer. Math. Soc.}, 136:825--837, 2008.

\bibitem{lencfm}
C.~Lenart.
\newblock On combinatorial formulas for {M}acdonald polynomials.
\newblock {\em Adv. Math.}, 220:324--340, 2009.

\bibitem{lenhhl}
C.~Lenart.
\newblock Haglund-{H}aiman-{L}oehr type formulas for {H}all-{L}ittlewood
  polynomials of type {$B$} and {$C$}.
\newblock {\em Algebra Number Theory}, 4:887--917, 2010.

\bibitem{lenhlp}
C.~Lenart.
\newblock Hall-{L}ittlewood polynomials, alcove walks, and fillings of {Y}oung
  diagrams ({A}ppendix with {A}. {L}ubovsky).
\newblock {\em Discrete Math.}, 311:258--275, 2011.

\bibitem{lapawg}
C.~Lenart and A.~Postnikov.
\newblock Affine {W}eyl groups in {$K$}-theory and representation theory.
\newblock {\em Int. Math. Res. Not.}, pages 1--65, 2007.
\newblock Art. ID rnm038.

\bibitem{lapcmc}
C.~Lenart and A.~Postnikov.
\newblock A combinatorial model for crystals of {K}ac-{M}oody algebras.
\newblock {\em Trans. Amer. Math. Soc.}, 360:4349--4381, 2008.

\bibitem{macsfh}
I.~G. Macdonald.
\newblock {\em Symmetric Functions and \mbox{H}all Polynomials}.
\newblock Oxford Mathematical Monographs. Oxford University Press, Oxford,
  second edition, 1995.

\bibitem{prorsl}
R.~Proctor.
\newblock Representations of {$\mathfrak{sl}(2,\,\mathbb{C})$} on posets and
  the {S}perner property.
\newblock {\em SIAM J. Algebraic Discrete Methods}, 3:275--280, 1982.

\bibitem{ramawh}
A.~Ram.
\newblock Alcove walks, {H}ecke algebras, spherical functions, crystals and
  column strict tableaux.
\newblock {\em Pure Appl. Math. Q.}, 2:963--1013, 2006.

\bibitem{raycfm}
A.~Ram and M.~Yip.
\newblock A combinatorial formula for {M}acdonald polynomials.
\newblock {\em Adv. Math.}, 226:309--331, 2011.

\bibitem{sasrsa}
B.~Sagan and R.~Stanley.
\newblock Robinson-{S}chensted algorithms for skew tableaux.
\newblock {\em J. Combin. Theory Ser. A}, 55:161--193, 1990.

\bibitem{schghl}
C.~Schwer.
\newblock Galleries, {H}all-{L}ittlewood polynomials, and structure constants
  of the spherical {H}ecke algebra.
\newblock {\em Int. Math. Res. Not.}, pages Art. ID 75395, 31, 2006.

\bibitem{stausa}
R.~Stanley.
\newblock Unimodal sequences arising from {L}ie algebras.
\newblock In {\em Combinatorics, representation theory and statistical methods
  in groups}, volume~57 of {\em Lecture Notes in Pure and Appl. Math.}, pages
  127--136. Dekker, New York, 1980.

\bibitem{stawgh}
R.~Stanley.
\newblock Weyl groups, the hard {L}efschetz theorem, and the {S}perner
  property.
\newblock {\em SIAM J. Algebraic Discrete Methods}, 1:168--184, 1980.

\bibitem{stagln}
R.~Stanley.
\newblock {${GL}(n,{\mathbb{C}})$} for combinatorialists.
\newblock In {\em Surveys in combinatorics ({S}outhampton, 1983)}, volume~82 of
  {\em London Math. Soc. Lecture Note Ser.}, pages 187--199. Cambridge Univ.
  Press, Cambridge, 1983.

\bibitem{stasbs}
R.~Stanley.
\newblock The stable behavior of some characters of {${SL}(n,\mathbb{C})$}.
\newblock {\em Linear and Multilinear Algebra}, 16:3--27, 1984.

\bibitem{stauls}
R.~Stanley.
\newblock Unimodality and {L}ie superalgebras.
\newblock {\em Stud. Appl. Math.}, 72:263--281, 1985.

\bibitem{stadp}
R.~Stanley.
\newblock Differential posets.
\newblock {\em J. Amer. Math. Soc.}, 1:919--961, 1988.

\bibitem{stalcu}
R.~Stanley.
\newblock Log-concave and unimodal sequences in algebra, combinatorics, and
  geometry.
\newblock In {\em Graph theory and its applications: {E}ast and {W}est
  ({J}inan, 1986)}, volume 576 of {\em Ann. New York Acad. Sci.}, pages
  500--535. New York Acad. Sci., New York, 1989.

\bibitem{stascp}
R.~Stanley.
\newblock Some combinatorial properties of {J}ack symmetric functions.
\newblock {\em Adv. Math.}, 77:76--115, 1989.

\bibitem{stavdp}
R.~Stanley.
\newblock Variations on differential posets.
\newblock In {\em Invariant theory and tableaux ({M}inneapolis, {MN}, 1988)},
  volume~19 of {\em IMA Vol. Math. Appl.}, pages 145--165. Springer, New York,
  1990.

\bibitem{staec2}
R.~Stanley.
\newblock {\em Enumerative {C}ombinatorics. {V}ol. 2}, volume~62 of {\em
  Cambridge Studies in Advanced Mathematics}.
\newblock Cambridge University Press, Cambridge, 1999.

\bibitem{stertt}
J.~Stembridge.
\newblock Rational tableaux and the tensor algebra of {$\mathfrak{gl}_n$}.
\newblock {\em J. Combin. Theory Ser. A}, 46:79--120, 1987.

\bibitem{leuajt}
M.~A.~A. van Leeuwen.
\newblock An analogue of jeu de taquin for {L}ittelmann's crystal paths.
\newblock {\em S\'em. Lothar. Combin.}, 41:Art.\ B41b, 23 pp.\ (electronic),
  1998.

\bibitem{yiplrr}
M.~Yip.
\newblock A {L}ittlewood-{R}ichardson rule for {M}acdonald polynomials.
\newblock {\em Math. Z.}, 272:1259--1290, 2012.

\end{thebibliography}


}
\end{document}